\documentclass[cmp,final]{svjour}
\usepackage{amsfonts}
\usepackage{hyperref}
\hypersetup{
colorlinks = true,
citecolor=blue,
urlcolor=blue,
linkcolor=blue,
pdfpagemode = UseNone
}

\newcommand{\lbl}[1]{{\tt \small [#1]}\label{#1}}
\renewcommand{\lbl}[1]{\label{#1}}

\newcommand{\rf}[1]{(\ref{#1})}
\newcommand{\rfg}{\langle f, g \rangle}
\newcommand{\rfh}{\langle f, h \rangle}
\newcommand{\rgh}{\langle g, h \rangle}

\newcommand{\aaa}{{\mathbf a}}
\newcommand{\XX}{{\mathbf X}}
\newcommand{\YY}{{\mathbf Y}}
\newcommand{\ZZ}{{\mathbf Z}}
\newcommand{\II}{{\mathbf I}}
\newcommand{\EE}{{\mathbb E}}

\newcommand{\sR}{{\mathbb R}}
\newcommand{\sC}{{\mathbb C}}

\newcommand{\vacuum}{\Phi}

\newcommand{\HH}{\Gamma_q({\mathcal H})} 
\renewcommand{\H}{{\mathcal H}} 
\newcommand{\X}{\tilde{X}}
\newcommand{\Y}{\tilde{Y}}
\newcommand{\Z}{\tilde{Z}}

\newcommand{\calA}{{\mathcal A}}
\newcommand{\calB}{{\mathcal B}}
\newcommand{\CC}{{\mathbf C}} 

\newenvironment{proofof}[1]{\noindent {\bf #1.\/}}{\qed\vskip 0.1in}

\author{
Wlodzimierz  Bryc
}
\institute{
Department of Mathematics,
University of Cincinnati,
PO Box 210025,
Cincinnati, OH 45221--0025, USA.
Wlodzimierz.Bryc@UC.edu
}
\date{
Received: 18 September, 2000/Accepted: 17 November, 2000
}
\communicated{H. Araki}

\title{Classical Versions of $q$-{G}aussian Processes: Conditional Moments and Bell's Inequality}
\titlerunning{$q$-{G}aussian Processes: Conditional Moments and Bell's Inequality}
 \usepackage{framed,color,pgf}
\begin{document}

\maketitle
\begin{abstract}
We show that classical processes corresponding to operators which satisfy a $q$-commutative
relation have linear regressions and quadratic conditional variances.
From this we deduce that Bell's inequality for their covariances can be extended from $q=-1$
to the entire range $-1\leq q <1$.

\end{abstract}
\colorlet{shadecolor}{gray!10}

\section*{Corrections of Sunday, February 29, 2004 at 15:40}
The following corrections were found  after the printed version  appeared in Comm. Math. Phys., 219 (2001), pp. 259-270.
\begin{enumerate}
\item Of course, the expected value $\EE:\calA\to\sC$, not into $\sR$.
\item Conditional variance in Proposition \ref{WWW} is now correct.
\end{enumerate}
 \begin{leftbar} \begin{shaded}
\section*{Corrections of Sunday, April 14, 2024  09:38}
\begin{enumerate}
\item Another typo in the conditional variance in Proposition \ref{WWW} was corrected.
\end{enumerate}
\end{shaded}\end{leftbar}

\section{Introduction}

In this paper we  consider a linear mapping $\H \ni f\mapsto \aaa_f\in\calB$ from the real Hilbert space $\H$
into the algebra $\calB$ of bounded operators acting on a complex Hilbert space 
which satisfies the  $q$-commutation relations
\begin{equation}\lbl{q-comute}
\aaa_f \aaa_g^*-q\aaa_g^*\aaa_f=\langle f,g\rangle \II,
\end{equation}
and $\aaa_f\vacuum=0$ for a vacuum vector $\vacuum$.
This defines a non-commutative stochastic process
$\XX_f=\aaa_f+\aaa_f^*$, first studied in \cite{FB70}, which following  \cite{BKS97} we
 call the $q$-Gaussian process. For different values of $q$, these
 processes  interpolate between
the bosonic ($q=1$) and fermionic ($q=-1$) processes, and include free processes of Voiculescu
\cite{VDN92}
($q=0$).

One of the basic problems arising in this context is the existence of the classical versions
of $q$-Gaussian processes, see Definition \ref{CV-def}. For $q=1$, these are the classical
Gaussian processes with the covariances $\rfg_{f,g\in\H}$.
For $q=-1$, the classical versions are two-valued, so
 Bell's inequality \cite{Bell64} shows that only some covariances may have the
classical versions. In \cite{FB70}  classical versions were constructed for
covariances corresponding to stationary two-valued
Markov processes ($q=-1$).
%
In \cite[Prop. 3.9]{BKS97}, the existence of such classical
versions was  proved for all $-1<q<1$ in the case where the $q$-Gaussian
process is Markovian (which can be characterized in
terms of the covariance function).

The situation
for other covariances remained open in \cite{BKS97} and
it was unclear which $q$-Gaussian processes have no classical
realizations.
This issue is addressed in the present paper. Using a formula for
 conditional variances
of classical versions  we derive a
 constraint on the covariance which extends
one of the Bell's inequalities from  $q=-1$ to general $-1\leq q<1$.
 The inequality
 implies that there are
covariances such that the corresponding
 non-commutative $q$-Gaussian processes cannot have  classical
versions over the entire range $-1\leq q<1$.
Since $q$ interpolates between the values $q=-1$,
 where classical versions may
 fail to exist and $q=1$,
where the classical versions always exist, it is
interesting that there is a version of Bell's inequality which  does not
depend on $q$.

 The proof relies on  formulas for conditional moments of the first two orders,
 which are of independent
interest.  Computations to derive them  were   possible thanks to
recent advances in  the  Fock space representation of   $q$-commutation
relations  \rf{q-comute}, see \cite{BKS97,BS91}.

\section{Preliminaries}
This section introduces the Fock space representation of $q$-Gaussian processes,
and states known results
in the  form convenient for us.
It is based on \cite{BKS97}.
\subsection{Notation}
Throughout the paper, $q$ is a fixed parameter and $-1<q<1$.
For $n=0, 1, 2,\dots$ we define $q$-integers
$
[n]_q:=\frac{1-q^n}{1-q} .
$
The $q$-factorials are
$
[n]_q!:=[1]_q[2]_q\dots[n]_q,
$
with the  convention $[0]_q!:=1$.

The $q$-Hermite polynomials are defined by the recurrence
\begin{equation}
\lbl{q-H}
xH_n(x)=H_{n+1}(x)+[n]_qH_{n-1}(x), \, n\geq 0
\end{equation}
with $H_{-1}(x):=0, H_0(x):=1$.
These polynomials
are orthogonal with respect to the unique absolutely continuous  probability measure
$\nu_q(dx)=f_q(x) dx$ supported on $[{-2/\sqrt{1-q}},{2/\sqrt{1-q}}]$,
where density
$f_q(x)$ has explicit product expansion, see \cite[Theorem 1.10]{BKS97} or \cite{KS98};
the  second moments of $q$-Hermite polynomials are
$\int_{-2/\sqrt{1-q}}^{2/\sqrt{1-q}} \left(H_n(x)\right)^2\nu_q(dx)=[n]_q!
$.
In our notation  we are suppressing the dependence of $H_n(x)$
on $q$.

\subsection{$q$-Fock space}
%

For a real Hilbert space $\H$ with complexification $\H_c=\H\oplus i\H$ we define its $q$-Fock space $\HH$ as the closure
 of $\sC\vacuum\oplus \bigoplus_n\H_c^{\otimes n}$, the
linear span of vectors $f_1\otimes\dots\otimes f_n$,
in the scalar product
\begin{equation} \lbl{q-scalar}
\langle f_1\otimes\dots\otimes  f_n|g_1\otimes\dots\otimes g_m\rangle_q =
\left\{\begin{array}{ll}
\sum_{\sigma\in S_n}q^{|\sigma|}\prod_{j=1}^n \langle f_j,g_{\sigma(j)}\rangle & \mbox{ if } m=n \\
0 &\mbox{ if } m\ne n
\end{array}
\right..
\end{equation}
Here $\vacuum$ is the vacuum vector,
 $S_n$ are permutations of $\{1,\dots,n\}$ and $|\sigma|=\#\{(i,j):i<j,
\sigma(i)>\sigma(j)\}$.
For the proof that \rf{q-scalar} indeed is non-negative definite,
 see \cite{BS91}.

Given the $q$-Fock space $\HH$  and $f\in\H$ we define the creation operator
$\aaa_f:\HH\to\HH$ and its $\langle\cdot|\cdot\rangle_q$-adjoint, the annihilation operator
 $\aaa_f^*:\HH\to\HH$ as follows:
$$
\aaa_f\vacuum:=0,
$$
\begin{equation}\lbl{anihilate}
\aaa_f f_1\otimes\dots\otimes
f_n:=\sum_{j=1}^n q^{j-1}\langle f, f_j \rangle f_1\otimes\dots\otimes f_{j-1}\otimes f_{j+1}\otimes\dots\otimes f_n,
\end{equation}
and
$$\aaa_f^*\vacuum=f,$$
\begin{equation}\lbl{create}
\aaa_f^* f_1\otimes\dots\otimes f_n:=f\otimes f_1\otimes\dots\otimes f_n.
\end{equation}
These operators are bounded,
satisfy commutation relation \rf{q-comute}, and  $\aaa_{f+g}=\aaa_f+\aaa_g$,
see \cite{BS91}.

\subsection{$q$-Gaussian processes}

We now consider (non-commutative) random variables as the elements of the algebra $\calA$
 generated
by the self-adjoint operators
$ \XX_f:=\aaa_f+\aaa_f^*
$, with vacuum expectation state $\EE:\calA\to\sC$ given by $\EE(\XX)=\langle \vacuum|\XX\vacuum \rangle_q$.

\begin{definition}
We will call $\{\XX({t}): {t\in T}\}$ a  $q$-Gaussian (non-commutative) process indexed by $T$
if
 there are vectors $h(t)\in \H$ such that $\XX(t)=\XX_{h(t)}$.
\end{definition}
For a $q$-Gaussian process  the covariance function
$c_{t,s}:=\EE(\XX_t\XX_s)$ becomes
$c_{t,s}=\langle h(t), h(s) \rangle $.

The Wick products $\psi(f_1\otimes\dots\otimes f_n)\in\calA$ are defined recurrently by
  $\psi(\vacuum):=\II$, and
\begin{eqnarray}\lbl{Wick-def}
&\psi(f\otimes f_1\otimes\dots\otimes f_n):=& \\ \nonumber
&\XX_f \psi(f_1\otimes\dots\otimes f_n) - \sum_{j=1}^n
q^{j-1}\langle f, f_j \rangle  \psi(f_1\otimes\dots\otimes f_{j-1}\otimes f_{j+1}\otimes\dots\otimes f_n).
&
\end{eqnarray}
An important property of Wick products is that if $\XX=\psi(f_1\otimes\dots\otimes f_n)$ then
\begin{equation}\lbl{Wick-appl}
\XX\vacuum=f_1\otimes\dots\otimes f_n.
\end{equation}
 We will also use the connection with $q$-Hermite polynomials.
If $\|f\|=1$
then
\begin{equation}\lbl{psi-Hermit}
\psi\left(f^{\otimes n}\right)=H_n\left(\XX_f\right),
\end{equation}
see \cite[Prop. 2.9]{BKS97}. Formulas \rf{q-scalar}, \rf{Wick-appl}, and  \rf{psi-Hermit}
show that for a unit vector $f\in\H$ we have
\begin{equation}\lbl{H-norm nc}
\EE\left(\left(H_n(\XX_f)\right)^2\right)=\sum_{\sigma\in S_n}q^{|\sigma|}=[n]_q!.
\end{equation} Thus
 $\nu_q$ is indeed the distribution of $\XX_f$.
Our main use of the Wick product is to compute certain conditional expectations.
\subsection{Conditional expectations}

Recall that a
(non-commutative) conditional expectation on the probability space $(\calA,\EE)$ with respect to the subalgebra $\cal B\subset \calA$
is a mapping $\cal E: \calA\to\cal B$ such that
\begin{equation}\lbl{E nc}
\EE(\YY_1 \XX \YY_2)=\EE(\YY_1{\cal E}(\XX) \YY_2)
\end{equation}
 for all $\XX\in\calA, \YY_1,\YY_2\in\cal B$.

We will study only  algebras $\cal B$ generated
 by the identity and the finite number
of random variables $\XX_{f_1},\dots,\XX_{f_n}$.
In this situation, we will use  a more probabilistic
notation:
$$\EE(\XX|\XX_{f_1},\dots,\XX_{f_n}):=\cal E(\XX), \,\XX\in\calA.$$

In this setting conditional expectations are easily computed for $\XX$ given by Wick products.
This important  result comes from \cite[Theorem 2.13]{BKS97}.
\begin{theorem}\lbl{psi-conditionals} If $\YY=\psi(g_1\otimes\dots \otimes g_m)$,
$\XX_1=\XX_{f_1},\dots,\XX_k=\XX_{f_k}$ for some $f_i,g_j\in\H$ and
$P:\H\to \H$ denotes orthogonal projection onto the span of $f_1,\dots,f_k$ then
$$
\EE(\YY|\XX_1,\dots,\XX_k)=\psi(Pg_1\otimes\dots\otimes P g_m).
$$
\end{theorem}

The following formula is an immediate consequence
of Theorem \ref{psi-conditionals} and \rf{psi-Hermit}, and is implicit in \cite[Proof
of Theorem 4.6]{BKS97}.
\begin{corollary}\lbl{eigen vectors}
If $\XX=\XX_f, \YY=\XX_g$ with unit vectors $\|f\|=\|g\|=1$ and $H_n$ is the $n^{\mbox{\small th}}$ $q$-Hermite polynomial, see \rf{q-H}, then
\begin{equation}\lbl{E(H|X)}
\EE(H_n(\YY)|\XX)=\rfg^n H_n(\XX).
\end{equation}
\end{corollary}

 For a  finite number of
vectors
$f_0,f_1,\dots,f_k\in\H$, let $\XX_k:=\XX_{f_k}$.
These
(non-commutative) random variables have
 linear regressions and constant conditional variances like
the classical (commutative) Gaussian random variables.

\begin{proposition}\lbl{q-regression}
\begin{equation}\lbl{LR}
\EE(\XX_0|\XX_1,\dots,\XX_k)=\sum_{j=1}^k a_j \XX_j
\end{equation}
and
\begin{equation}\lbl{QV}
\EE(\XX_0^2|\XX_1,\dots,\XX_k)=(\sum_{j=1}^k a_j \XX_j)^2+c\II.
\end{equation}
If $f_1\dots,f_k\in\H$ are linearly independent then the
  coefficients $a_j, c$ are uniquely determined by the covariance matrix
$C=[c_{i,j}]:=[\langle f_i,f_j\rangle]$.
\end{proposition}
Notice that Eq. \rf{QV} can indeed
be rewritten as the statement that conditional variance is constant,

$$Var(\XX_0|\XX_1,\dots,\XX_k):=\EE\left(\left(\XX_0-\EE(\XX_0|\XX_1,\dots,\XX_k)\right)^2|\XX_1,\dots,\XX_k\right)=c\II.$$
\begin{proof}
This follows from Theorem \ref{psi-conditionals} and \rf{Wick-def}.
Write the orthogonal projection of $f_0$ onto the span of $f_1,\dots f_k$
as the linear combination $g=\sum_ja_j f_j$. Then
$\EE(\XX_0|\XX_1,\dots,\XX_k)=E(\psi(f_0)|\XX_1,\dots,\XX_k)=\psi(g)=\sum_j a_j\psi(f_j)$,
which proves \rf{LR}. Similarly,
$
\EE(\XX_0^2-\|f_0\|^2\II|\XX_1,\dots,\XX_k)=\EE(\psi(f_0\otimes f_0)|\XX_1,\dots,\XX_k)=
\psi(g\otimes g)=(\sum_j a_j\XX_j)^2-\|g\|^2\II$.
This proves \rf{QV} with $c=\|f_0\|^2-\|g\|^2$.

If $f_1\dots,f_k\in\H$ are linearly independent then the representation
 $g=\sum_ja_j f_j$ is unique.
\end{proof}
To   analyze standardized triplets in more detail
 we   need the explicit form of the coefficients. (We omit the straightforward
calculation.)
\begin{corollary}\lbl{q-reg3}
If $\XX:=\XX_f, \YY:=\XX_g,\ZZ:=\XX_h$ and $f,h\in\H$ are linearly independent unit vectors, then
\begin{equation}\lbl{LR3}
\EE(\YY|\XX,\ZZ)=a\XX+b\ZZ,
\end{equation}
\begin{equation}\lbl{QV3}
\EE(\YY^2|\XX,\ZZ)=\left(a\XX+b\ZZ\right)^2+c\II,
\end{equation}
where
\begin{eqnarray}
\lbl{a} a&=& \frac{\langle f, g \rangle - \langle g , h\rangle \langle f, h \rangle}{
 1-\langle f, h \rangle^2}, \\
\lbl{b} b&=&\frac{\langle g, h \rangle - \langle f , g\rangle \langle f, h\rangle}{
 1-\langle f, h \rangle^2}.
\end{eqnarray}
\end{corollary}

Another calculation shows that $c=\det(\CC)/( 1-\langle f, h \rangle^2)$, where $\CC$ is the
covariance matrix; in particular
$c\geq
0$.
\section{Conditional Moments of Classical Versions}

We give the definition of a classical version
which is convenient for bounded processes; for a more general definition, see
\cite[Def. 3.1]{BKS97}.

\begin{definition}\lbl{CV-def} A classical version of the process $\XX(t)$ indexed by
${t\in T\subset \sR}$
is a stochastic process $\X(t)$ defined on some classical probability space such that for
any finite number of indexes $t_1<t_2<\dots<t_k$ and any polynomials
$P_1,\dots,P_k$,
\begin{eqnarray}\lbl{C-ver}
\EE\left(P_1(\XX(t_1))P_2(\XX(t_2))\dots P_k(\XX(t_k))\right)=
\\ \nonumber
E\left(P_1(\X(t_1))P_2(\X(t_2))\dots P_k(\X(t_k))\right).
\end{eqnarray}
\end{definition}
Here $E(\cdot)$ denotes the classical expected value given by Lebesgue integral with respect to
the classical probability measure.

Our main interest is in finite index set $T=\{t_1,t_2,t_3\}$, where $t_1<t_2<t_3$.
In this case we  write
$\XX:=\XX(t_1), \YY:=\XX(t_2),\ZZ:=\XX(t_3)$. We
 say that an ordered triplet $(\XX,\YY,\ZZ)$ has a classical version $\X,\Y,\Z$,
if $\EE\left(P_1(\XX)P_2(\YY)P_3(\ZZ)\right)=E\left(P_1(\X)P_2(\Y)P_3(\Z)\right)$
for all polynomials $P_1,P_2,P_3$.

The classical version of a non-commutative process is order-dependent,
 since the left-hand side of \rf{C-ver} may depend on
the ordering of the variables, while the right-hand side does not.
 For specific example
 in the context of
$q$-Gaussian random variables,
 see \cite[formulas (2,64) and (2.65)]{FB70}.

\subsection{Triplets}
All pairs $(\XX_f,\XX_g)$ of $q$-Gaussian random variables have classical versions
because $\EE(\XX_f^m\XX_g^n)=\EE(\XX_f^n\XX_g^m)$ for all integer $m,n$;
however, the classical version of a triplet  may fail to exist.
With this in mind we consider $q$-Gaussian triplets
\begin{equation}
\lbl{triplet}
\XX:=\XX_f, \YY:=\XX_g,\ZZ:=\XX_h.
\end{equation}
To  simplify the notation we take unit vectors
$\|f\|=\|g\|=\|h\|=1$. We   assume that there is
a classical version $(\X,\Y,\Z)$ of $(\XX,\YY,\ZZ)$, in this order.

From Corollary \ref{q-reg3} we know that
non-commutative random variables $\XX,\YY,\ZZ$ have linear regression and
constant conditional variance.
It turns out that the corresponding
classical random variables $\X,\Y,\Z$ also
have linear regressions, while
their conditional variances get  perturbed into quadratic polynomials.
\begin{theorem}\lbl{q-reg classic}
If $(\X,\Y,\Z)$ is a classical version of the $q$-Gaussian triplet \rf{triplet} then
\begin{equation}
\lbl{LR-c}
E(\Y|\X,\Z)=a\X+b\Z,\end{equation}
\begin{equation}\lbl{QV-c}
E(\Y^2|\X,\Z)=A\X^2+B\X\Z+C\Z^2+D,
\end{equation}
where
 $a,b$ are given by \rf{a},\rf{b},
\begin{eqnarray}
\lbl{A} A&=&\frac{ab(1-q)\rfh+a^2(1-q\rfh^2)}{1-q\rfh^2},
\\
\lbl{B} B&=&\frac{ab(1+q)(1-\rfh^2)}{1-q\rfh^2},
\\
\lbl{C} C&=&\frac{ab(1-q)\rfh+b^2(1-q\rfh^2)}{1-q\rfh^2},
\end{eqnarray}
 and
\begin{equation}\lbl{D}
D=1-A-B\rfh -C.
\end{equation}
\end{theorem}


The proof relies on the following technical result.
\begin{lemma}\lbl{switch}
If $H_n,H_m$ are $q$-Hermite polynomials given by \rf{q-H},
then

\begin{eqnarray}
\lbl{H(X)ZXH(Z)} \\ \nonumber
\EE\left(H_n(\XX)\ZZ\XX H_m(\ZZ)\right)&=&
\left\{\begin{array}{ll}
\rfh^{n+1}[n+2]_q! & \mbox{ if } m=n+2 \\
\rfh^{n-1}[n]_q! & \mbox{ if } m=n-2 \\
\rfh^{n-1}\left(([n]_q+1)\rfh^2+q[n]_q\right)[n]_q! & \mbox{ if } m=n \\
0 & \mbox{ otherwise }
\end{array}\right.,
\\
\lbl{H(X)XZH(Z)} \\ \nonumber
\EE\left(H_n(\XX)\XX\ZZ H_{m}(\ZZ)\right)&=&
\left\{\begin{array}{ll}
\rfh^{n+1}[n+2]_q! & \mbox{ if } m=n+2 \\
\rfh^{n-1}[n]_q! & \mbox{ if } m=n-2 \\
\rfh^{n-1}\left([n+1]_q\rfh^2+[n]_q\right)[n]_q! & \mbox{ if } m=n \\
0 & \mbox{ otherwise }
\end{array}\right. ,
\\\label{H(X)Z2H(Z)}\lbl{H(X)X2H(Z)}
\EE\left(H_n(\XX)\XX^2H_m(\ZZ)\right)&=&
\EE\left(H_m(\XX)\ZZ^2H_n(\ZZ)\right)= \\ \nonumber &&
\left\{\begin{array}{ll}
\rfh^{n+2}[n+2]_q! & \mbox{ if } m=n+2 \\
\rfh^{n-2}[n]_q! & \mbox{ if } m=n-2 \\
{\rfh^n\left([n+1]_q+[n]_q\right)[n]_q!} & \mbox{ if } m=n \\
0 & \mbox{ otherwise }
\end{array}\right. ,
\end{eqnarray}
\begin{eqnarray}\lbl{L2}\EE\left(H_n(\XX) H_m(\ZZ)\right)&=&
\left\{\begin{array}{ll}
\rfh^n[n]_q! & \mbox{ if } m=n \\
0 & \mbox{ otherwise }
\end{array}\right. .
\end{eqnarray}
\end{lemma}
\begin{proofof}{\it Proof of Theorem \protect{\ref{q-reg classic}}}
Since $\X,\Y,\Z$ are bounded random variables, to prove \rf{LR-c} we need only to verify
that for arbitrary polynomials $P,Q$,
$$
E\left(P(\X)\Y Q(\Z)\right)=E\left(P (\X)(a\X+b\Z) Q (\Z)\right).$$

 This is equivalent to
$$\EE\left(P(\XX)\YY Q(\ZZ)\right)=\EE\left(P (\XX)(a\XX+b\ZZ) Q (\ZZ)\right),
$$
see  \rf{C-ver}.
The latter follows from \rf{LR3} and \rf{E nc},  proving \rf{LR-c}.

To prove \rf{QV-c}, we  verify
that for arbitrary polynomials $P,Q$ we have

$$E\left(P(\X)\Y^2 Q(\Z)\right)=E\left(P (\X)(A\X^2+B\X\Z+C\Z^2+D) Q (\Z)\right).$$

By definition  \rf{C-ver}, this is equivalent to

\begin{equation}\lbl{QV-nc}
\EE\left(P(\XX)\YY^2 Q(\ZZ)\right)=\EE\left(P (\XX)(A\XX^2+B\XX\ZZ+C\ZZ^2+D) Q (\ZZ)\right).
\end{equation}
It suffices to show that \rf{QV-nc} holds true when
$P=H_n$ and $Q=H_m$ are the $q$-Hermite polynomials defined by \rf{q-H}.
Formula \rf{QV3} implies that the left-hand side of \rf{QV-nc} is given by
$$ c\EE(H_n(\XX)H_m(\ZZ))+
a^2\EE(H_n(\XX)\XX^2H_m(\ZZ))+b^2 \EE(H_n(\XX)\ZZ^2H_m(\ZZ))$$
$$
+ab \EE(H_n(\XX)\XX\ZZ H_m(\ZZ))+ab \EE(H_n(\XX)\ZZ\XX H_m(\ZZ)),
$$
and the right-hand side becomes
$$
A\EE(H_n(\XX)\XX^2H_m(\ZZ))+C \EE(H_n(\XX)\ZZ^2H_m(\ZZ))
+$$
$$B \EE(H_n(\XX)\XX\ZZ H_m(\ZZ))+ D\EE(H_n(\XX)H_m(\ZZ)).
$$
Using formulas from Lemma \ref{switch} we can see that both sides are zero, except when $m=n$
or $m=n\pm 2$. We now consider these three cases separately.

\noindent
{\it Case $m=n+2$}:
Using Lemma \ref{switch},  \rf{QV-nc} simplifies to
$$
(a^2\rfh^2+2ab \rfh+b^2)\rfh^n[n+2]_q!=(A\rfh^2+B\rfh+C)\rfh^n[n+2]_q!.
$$
This equation is satisfied when   coefficients $A,B,C$ satisfy the equation
\begin{equation}\lbl{ABC 1}
A\rfh^2+B\rfh+C=a^2\rfh^2+2ab \rfh+b^2.
\end{equation}

\noindent
{\it Case $m=n-2$}:
Using Lemma \ref{switch},  \rf{QV-nc} simplifies to
$$
(a^2+2ab \rfh+b^2\rfh^2)\rfh^{n-2}[n]_q!=(A+B\rfh+C\rfh^2)\rfh^{n-2}[n]_q!.
$$
This equation is satisfied whenever
\begin{equation}\lbl{ABC 2}
A+B\rfh+C\rfh^2=a^2+2ab \rfh+b^2\rfh^2.
\end{equation}

\noindent {\it Case $m=n$}: 
We use again Lemma \ref{switch}.
On both sides of Eq. \rf{QV-nc} we factor out
$\rfh^{n-1}[n]_q!$, and equate the remaining coefficients.
(This is allowed since we are after sufficient conditions only!)
We get
$$
\left(\rfh (a^2+b^2)([n+1]_q+[n]_q)+ab \left(\rfh^2[n+1]_q+
(1+q)[n]_q\right)+c\rfh
\right)
=$$
$$
\left((1+q)(A+C)+B(q\rfh^2+1)[n]_q+D\rfh\right).
$$
Now we use $[n+1]_q=1+q[n]_q$. Suppressing the correction to
the constant term (i.e., the term free of $n$),
we get

$$
(1+q)\left(\rfh (a^2+b^2)+ab(1+\rfh^2)\right)[n]_q+c\rfh+\dots
=$$

$$
\left((1+q)(A+C)\rfh+B(q\rfh^2+1)\right)
[n]_q+D\rfh,
$$
where $c+\dots$  denotes the suppressed constant term corrections.

This equation holds true when  the coefficients at $[n]_q$ match, which gives
\begin{equation}\lbl{ABC 3}
(1+q)\rfh(A+C)+B(q\rfh^2+1)=(1+q)\left((a^2+b^2)\rfh+ab(\rfh^2+1)\right),
\end{equation}
and the constant terms match: $c+\dots=D$. The latter holds true
 when the expectations are equal  ($n=m=0$), and hence this condition
is equivalent to
\rf{D}. The remaining three equations  \rf{ABC 1}, \rf{ABC 2}, and \rf{ABC 3} have
a unique solution given by the expressions \rf{A}, \rf{B}, \rf{C}.
\end{proofof}

\begin{proofof}{\it Proof of Lemma \protect{\ref{switch}}}
 Using the definition of vacuum expectation state,
\rf{Wick-appl} and \rf{psi-Hermit} we get
$\EE(H_n(\XX)\ZZ\XX H_m(\ZZ))=
\langle \ZZ H_n(\XX)\vacuum |\XX H_m(\ZZ)\vacuum \rangle_q=
\langle \XX_h f^{\otimes n}|\XX_f h^{\otimes m}\rangle_q$.
Therefore \rf{anihilate}, and \rf{create}   imply

\begin{eqnarray}\lbl{foo*}
&&\EE(H_n(\XX)\ZZ\XX H_m(\ZZ))= \\
&&\langle [n]_q\rfh f^{\otimes n-1}+h\otimes  f^{\otimes n}|
 [m]_q\rfh h^{\otimes m-1}+f\otimes  h^{\otimes m}
\rangle_q .\nonumber
\end{eqnarray}
The latter is zero, except when $m=n$ or $m=n\pm 2$.
We will consider these two cases separately.

If $m=n$, by orthogonality  we have
$$\EE(H_n(\XX)\ZZ\XX H_m(\ZZ))=
[n]_q^2\rfh^2
\langle f^{\otimes n-1}| h^{\otimes n-1}
\rangle_q
+
\langle h\otimes  f^{\otimes n}|f\otimes  h^{\otimes n}
\rangle_q.
$$
Clearly, $\langle f^{\otimes n-1}| h^{\otimes n-1}
\rangle_q =\rfh^{n-1}[n-1]_q!$; this can be seen either from \rf{H-norm nc} and \rf{E(H|X)},
or directly from the definition \rf{q-scalar}.

By \rf{q-scalar} the second term splits into the sum over permutations
$\sigma'\in S_{n+1}$ such that $\sigma'(1)=1$ and the sum over the permutations such that
$\sigma'(1)=k>1$. This gives

$$\langle h\otimes  f^{\otimes n}|f\otimes  h^{\otimes n}
\rangle_q=
\sum_{\sigma\in S_n} \rfh q^{|\sigma|}\rfh^n
+\sum_{k=2}^{n+1}\sum_{\sigma\in S_n}q^{k-1+|\sigma|}\rfh^{n-1}
=$$
$$\rfh^{n+1}[n]_q!+\rfh^{n-1}q[n]_q [n]_q!.
$$
Elementary algebra now yields \rf{H(X)ZXH(Z)} for $m=n$.

If $m=n+2$, then the right-hand side of \rf{foo*} consists of only one term
we get
$$\EE(H_n(\XX)\ZZ\XX H_{n+2}(\ZZ))=
 [n+2]_q\rfh \langle h\otimes f^{\otimes n}| h^{\otimes n+1}
\rangle_q=$$
$$ [n+2]_q\rfh \sum_{\sigma\in S_{n+1}} q^{|\sigma|}\rfh^n=\rfh^{n+1}[n+2]_q!.
$$

Since $m=n-2$ is given by the same expression with the roles of $m,n$ switched around, this ends the proof of \rf{H(X)ZXH(Z)}.

The remaining expectations match the corresponding commutative values,
 and can also be  evaluated using  recurrence \rf{q-H} and formulas \rf{E(H|X)}, and
\rf{H-norm nc}.

To prove \rf{H(X)XZH(Z)} notice that
since $\XX$ and $H_n(\XX)$ commute, using \rf{q-H} and \rf{E(H|X)} we get
$$\EE(H_n(\XX)\XX\ZZ H_m(\ZZ))=\EE(\XX H_n(\XX)(H_{m+1}(\ZZ)+[m]_qH_{m-1}(\ZZ)))
=$$
$$
\rfh^{m+1}\EE(\XX H_n(\XX)H_{m+1}(\XX))+ [m]_q\rfh^{m-1}\EE(\XX H_n(\XX)H_{m-1}(\XX)).
$$
The only non-zero values are when $m=n$, or $m=n\pm 2$. Using \rf{q-H} again, and then
 \rf{H-norm nc}
we get \rf{H(X)XZH(Z)}.


%
Since  by \rf{E(H|X)} we have
$$\EE(H_n(\XX)\XX^2 H_m(\ZZ))=\rfh^m \EE\XX^2 H_n(\XX)H_m(\XX),$$
recurrence \rf{q-H} used twice proves \rf{H(X)X2H(Z)}.

Formula \rf{L2} is an immediate consequence of \rf{E(H|X)} and \rf{H-norm nc}.
\end{proofof}

\subsection{Relation to processes with independent increments}

In \cite[Definition 3.5]{BKS97} the authors define the non-commutative
$q$-Brownian motion and show that it has a classical version, see
\cite[Cor. 4.5]{BKS97}.
Since the classical version of the $q$-Brownian motion is Markov, Theorem \ref{q-reg classic} implies
that all regressions are linear, and all conditional variances are quadratic.
A   computation gives the following expression
for the conditional variances.


\begin{proposition}\lbl{WWW}
Let $\X_t$ be the classical version of the $q$-Brownian motion, i.~e.,
$\langle f_t,f_s\rangle =\min\{s,t\}$.
Then   for   $t_1<t_2<\dots<t_{n}<s<t$ we have
$$Var(\X_{s}|\X_{t_1},\dots,\X_{t_n},X_t)=$$
$${\frac {\left ({ t}-{ s}
\right )\left ({ s}-{ t_n}\right )}{\left ({ t}-q{ t_n}
\right )}}\left(1+{\frac { \left ({ \X_{t}}-{ \X_{t_n}
}\right )\left
({ t_n} { \X_{t}-{ t} { \X_{t_n}}}
\right )\left (1-q\right )}{\left ({
 t}-{ t_n}\right )^{2}}}\right).
$$
\end{proposition}
 In \cite{Wesolowski93},   classical processes with
independent increments, linear regressions, and quadratic conditional
variances are analyzed.
These processes have the same covariances as $q$-Brownian motion, but the conditional variances
are quadratic functions of the increment $\X_t-\X_{t_n}$ only.
Proposition \ref{WWW} shows that the classical realizations of $q$-Brownian motion are not among the processes
in \cite{Wesolowski93} and thus have dependent increments.

\section{Bell's Inequality}
It is   well known that all $q$-Gaussian $n$-tuples with $q=1$ have classical versions: these
are given by the classical Gaussian distribution with the same covariance matrix $[\langle
f_i,f_j\rangle]$.

For $q=-1$ the classical version of the the standardized
$q$-Gaussian
triplet $(\XX,\YY,\ZZ)$ consists of the $\pm 1$-valued symmetric random variables.
The celebrated Bell's inequality \cite{Bell64} therefore restricts their covariances:
\begin{equation}\lbl{Bell}
1-\langle f,h \rangle\geq |\langle f,g \rangle -\langle g,h \rangle|.
\end{equation}
In particular,  there are triplets of
$q$-Gaussian random variables with $q=-1$ which do not have a classical version.

The following shows that
 restriction \rf{Bell} is in force  for sub-Markov covariances over the entire range  $-1\leq q<1$.

\begin{theorem}\lbl{q-Bell}
Suppose that $(\X,\Y,\Z)$ is a classical version of $q$-Gaussian
$(\XX,\YY,\ZZ):=(\XX_f,\XX_g,\XX_h)$, where $f,h\in\H$ are linearly independent,
and $-1\leq q<1$.
If either \begin{equation}\lbl{sub-M}
\langle f, g\rangle \langle g, h \rangle \leq \langle f, h\rangle \mbox{ and }  0<\rfh<1,
\end{equation}
or $\langle f, h\rangle=0$, or $q=-1$, then inequality
\rf{Bell} holds true.
\end{theorem}
\begin{proof} Since the case $q=-1$ is well known, we restrict our attention to the case
 $-1<q<1$.
Our starting point is expression \rf{QV-c}.
A computation shows that
 the conditional variance $Var(\Y|\X,\Z):=E(\Y^2|\X,\Z)-\left(E(\Y|\X,\Z)\right)^2$ is as follows.
\begin{eqnarray} \lbl{foo1}
Var(\Y|\X,\Z)=&1-a^2-b^2-
ab\rfh \frac{(1+q)(1-\rfh^2)+2(1-q)}{1-q\rfh^2}
-&\\ \nonumber
&\frac {ab\left(1-q\right)}
{1-q\rfh^2}\left(\Z{\rfh-\X}\right )\left ({\rfh}\X-\Z\right).
&
\end{eqnarray}
The right-hand side of this expression must be non-negative over the support of $\X,\Z$.
It is known, see  \cite[Lemma
8.1]{Bryc98} or \cite[Theorem 1.10]{BKS97}, that $\X,\Z$ have the joint
probability density function $f(x,z)$ with respect to the product of marginals $\nu_q$.
Moreover, $f$ is
 defined for all $-2/\sqrt{1-q}\leq x,z\leq 2/\sqrt{1-q}$ and from
its explicit product expansion we can see that
 $$f(x,z)\geq \prod_{k=0}^\infty \frac{(1-\rfh^2 q^k)}{
(1+\rfh q^{k})^4}$$
 is  strictly positive. In particular, the right-hand side of \rf{foo1} must be non-negative
when evaluated at $\X=\sqrt{2}/\sqrt{1-q},\Z=-\sqrt{2}/\sqrt{1-q}$.

Using formulas \rf{a}, \rf{b} with the above values of $\X,\Y$
we get the rational expression for the conditional variance which can be written as follows.
\begin{eqnarray} \lbl{foo2}
&\left(1-q\rfh^2\right)\left(1-\rfh\right)^2Var\left(\Y|\X,\Z\right)=& \\ \nonumber
&
(1-\rfh)^2\left(1-q\rfh^2+(1+q)\rfh\rfg\rgh\right)-&
\\ \nonumber &\left(\rfg-\rgh\right)^2(1+\rfh^2) .
&
\end{eqnarray}
Therefore
$$(1-\rfh)^2\left(1-q\rfh^2+(1+q)\rfh\rfg\rgh\right)
\geq $$
$$\left(\rfg-\rgh\right)^2(1+\rfh^2).$$

Since the assumptions imply that $1-q\rfh^2+(1+q)\rfh\rfg\rgh\leq 1+\rfh^2$, this implies
$(1-\rfh)^2 \geq \left(\rfg-\rgh\right)^2$, proving \rf{Bell}.

\end{proof}

\subsection{Examples}
The first example shows that there are covariances such that $q$-Gaussian
random variables have no classical version for all $-1\leq q<1$.
\begin{example}\lbl{E1} Consider the case $\rfh=\rgh>0$. This can be realized when  the
covariance matrix is non-negative definite; a computation shows that this is equivalent
to the condition $2\rfh^2\leq 1+\rfg$. Since \rf{sub-M} is satisfied, Bell's inequality \rf{Bell}
implies $1+\rfg\geq 2\rfh$. Therefore, all choices of vectors $f,g,h\in\H$
such that
$
\rfh=\rgh$, $0<\rfh<1$, and $2\rfh^2-1<\rfg<2\rfh-1$
lead to $q$-Gaussian triplets with no classical version for $-1\leq q<1$.
\end{example}
A nice feature of Theorem \ref{q-Bell} is that its statement does not depend on $q$, as long as
$q<1$. But such a result cannot be sharp. A less transparent
statement
that the conditional variance is non-negative is a stronger restriction on the covariances,
and it depends on $q$. This is illustrated by the next example.
\begin{example}\lbl{E2}
Suppose $\rfh=\rgh=1/2$. Inequality \rf{Bell} used in Example \ref{E1}
implies that if a classical version of a $q$-Gaussian process exists
then $\rfg \geq 0$. Evaluating the conditional variance $Var(\Y|\X,\Z)$
at $\X=2/\sqrt{1-q},\Z=-\X$ we get a more restrictive constraint
$\rfg\geq \frac{q+5}{36}$.
\end{example}

\subsection*{Acknowledgements}
{\small I would like to thank the referee,  A.
Dembo, and T. Hodges for suggestions which improved the presentation, and to V. Kaftal
for several helpful discussions.
}


\end{document}